\documentclass{article}
\usepackage{verbatim}
\usepackage{amsmath}
\usepackage[mathscr]{eucal}
\usepackage{amsthm}
\usepackage{amssymb}
\usepackage[german,english]{babel}
\theoremstyle{plain}
\newtheorem{Theo}{Theorem}[section]
\newtheorem{Prop}[Theo]{Proposition}

\theoremstyle{definition}

\theoremstyle{remark}
\newtheorem{Rem}[Theo]{Remark}
\newcommand{\T}{\bar{T}}

\newcommand{\gl}{\mathfrak{g}\mathfrak{l}}
\newcommand{\g}{\mathfrak{g}}
\newcommand{\Par}{\operatorname{Par}}
\makeatletter
	
	\@addtoreset{equation}{section}
\makeatother
\begin{document}
\title{Schur-Weyl reciprocity between the quantum superalgebra 
and the Iwahori-Hecke algebra}
\author{Hideo Mitsuhashi}

\date{}

\maketitle

\begin{center}
Department of Information Technology \\
Kanagawa Prefectural Junior College for Industrial Technology \\
2--4--1 Nakao, Asahi--ku, Yokohama--shi, Kanagawa--ken 241--0815, Japan
\end{center}

\begin{abstract}
In this paper, we establish Schur-Weyl reciprocity between the quantum general super Lie algebra 
$U_q^\sigma\big{(}{\gl}(m,n)\big{)}$ and the Iwahori-Hecke algebra $\mathcal{H}_{\mathbb{Q}(q),r}(q)$. 
We introduce the sign $q$-permutation representation of $\mathcal{H}_{\mathbb{Q}(q),r}(q)$ 
on the tensor space $V^{{\otimes}r}$ of $(m+n)$ dimensional $\mathbb{Z}_2$-graded 
$\mathbb{Q}(q)$-vector space $V=V_{\bar{0}}{\oplus}V_{\bar{1}}$. 
This action commutes with that of $U_q^\sigma\big{(}{\gl}(m,n)\big{)}$ derived from 
the vector representation on $V$. 
Those two subalgebras of $\operatorname{End}_{\mathbb{Q}(q)}(V^{{\otimes}r})$ satisfy 
Schur-Weyl reciprocity. As special cases, we obtain the super case ($q{\rightarrow}1$), 
and the quantum case ($n=0$). Hence this result includes both the super case 
and the quantum case, and unifies those two important cases. 
\end{abstract}

\pagestyle{plain}
\pagenumbering{arabic}
\renewcommand{\theenumi}{\arabic{enumi}}
\renewcommand{\labelenumi}{(\theenumi)}

\section{Introduction}
In the representation theory, the classification and the construction of the irreducible 
representations are essential themes. 
In the first half of the twentieth century, I. Schur\cite{Schur} introduced a prominent method to 
obtain the finite dimensional irreducible representations of the general linear group 
$\operatorname{GL}(n,\mathbb{C})$, or equivalently of its Lie algebra $\gl(n,\mathbb{C})$, 
which we call Schur-Weyl reciprocity at present. 
Schur applied this method to the permutation action of the symmetric group $\mathfrak{S}_r$ 
and the diagonal action of $\operatorname{GL}(n,\mathbb{C})$ on the tensor powers $V^{{\otimes}r}$ 
of the $n$ dimensional complex vector space $V$. \par
After this work, Schur-Weyl reciprocity has been extended to various groups and algebras. 
Brauer\cite{Bra} obtained the centralizer algebra of the orthogonal Lie group $O(n)$. 
Sergeev\cite{Ser} and Berele-Regev\cite{B-R} extended the Schur's result to the general super 
Lie algebra $\gl(m,n)$. 
Jimbo\cite{Jimbo} extended it to the $q$-analogue case. He established Schur-Weyl reciprocity 
between the quantum enveloping algebra $U_q(\gl_{n+1})$ and the Iwahori-Hecke algebra of type $A$. 
As in the book of Curtis-Reiner\cite{C-R2}, the representation theory of Iwahori-Hecke algebras 
is an important part in representation theories of finite groups of Lie type. 
Hence we will focus on the representation theory of Iwahori-Hecke algebras. \par
In \cite{Mit}, we defined a $q$-deformation of the alternating group as a subalgebra of 
the Iwahori-Hecke algebra, and determined all the isomorphism classes of 
(ordinary) irreducible representations. 
After \cite{Mit}, we intended to compute character values of irreducible representations directly 
using combinatorial method. But this strategy did not go well because the notion of 
conjugacy classes of the $q$-deformation of the alternating group is obscure, hence we could not 
apply the classical($q=1$) case which is found in \cite{Headley}. 
Thereby we will take a representational approach to obtain the character table. \par
In \cite{Regev}, Regev obtained double centralizer properties for alternating groups. 
Roughly speaking, his result claims when one restrict the representation of the symmetric group 
on tensor space to the alternating group, 
the corresponding centralizer algebra enlarges in \lq\lq super case" while does not change in 
\lq\lq normal case".
Those facts suggest that Schur-Weyl reciprocity for $\gl(m,n)$ is more suitable to describe 
the representation theory of the alternating group than that for $\gl(n)$. \par
In this paper, we establish Schur-Weyl reciprocity between the quantum superalgebra 
$U_q^\sigma\big{(}{\gl}(m,n)\big{)}$ and the Iwahori-Hecke algebra $\mathcal{H}_{\mathbb{Q}(q),r}(q)$. 
We define the sign $q$-permutation representation of $\mathcal{H}_{\mathbb{Q}(q),r}(q)$ 
on $V^{{\otimes}r}$ using an operator $T$ on $V^{{\otimes}2}$ defined by: 
\begin{equation}
v_k{\otimes}v_lT=
\begin{cases}
\dfrac{(-1)^{|v_k|}(q+q^{-1})+q-q^{-1}}{2}v_k{\otimes}v_l& \text{if $k=l$,}\\
(-1)^{|v_k||v_l|}v_l{\otimes}v_k+(q-q^{-1})v_k{\otimes}v_l& \text{if $k<l$,}\\
(-1)^{|v_k||v_l|}v_l{\otimes}v_k& \text{if $k>l$,}
\end{cases}
\end{equation}
where $V$ is an $(m+n)$ dimensional $\mathbb{Z}_2$-graded $\mathbb{Q}(q)$-vector space and 
$|{\cdot}|$ is the degree map. 
This action reduces to the sign permutation action(see \cite{B-R}) of the symmetric group 
when $q{\rightarrow}1$ 
and to the well-known action of $\mathcal{H}_{\mathbb{Q}(q),r}(q)$ obtained from Drinfeld-Jimbo 
solutions to the Yang-Baxter equation when $n=0$. \par
The quantum superalgebra has been defined in several articles such as \cite{B-K-K},\cite{K-T} 
or \cite{Yamane}. 
$U_q^\sigma\big{(}{\gl}(m,n)\big{)}$ is a Hopf algebra obtained from the \lq\lq naive" quantum 
superalgebra $U_q\big{(}{\gl}(m,n)\big{)}$, which is a Hopf superalgebra, 
by adding an involutive element $\sigma$. 
We show that the vector representation of $U_q^\sigma\big{(}{\gl}(m,n)\big{)}$ on $V^{{\otimes}r}$, 
which is found in \cite{B-K-K}, and the sign $q$-permutation representation of 
$\mathcal{H}_{\mathbb{Q}(q),r}(q)$ are commutants of one another in $V^{{\otimes}r}$. 
Furthermore, extending the base field to the algebraic closure and applying 
the double centralizer theorem, we obtain the tensor space decomposition of 
${\big{(}V{\otimes}\overline{\mathbb{Q}(q)}}\big{)}^{{\otimes}r}$ as 
$\mathcal{H}_{\overline{\mathbb{Q}(q)},r}(q){\otimes}\big{(}U_q^\sigma\big{(}{\gl}(m,n)\big{)}
{\otimes}\overline{\mathbb{Q}(q)}\big{)}$-modules. 
Our result will be very useful to the representation theory of the $q$-deformation of the 
alternating group. 

\section{The sign $q$-permutation representation of the Iwahori-Hecke algebra of type $A$}
In this section, we shall define the sign $q$-permutation representation of 
the Iwahori-Hecke algebra of type $A$, which is a $q$-deformation of the sign permutation module 
introduced by several precedent works such as \cite{B-R}, \cite{Ser}. \par
Let $(W,S=\{s_1,\hdots,s_r\})$ be a Coxeter system of rank $r$. 
Let $R$ be a commutative domain with $1$, and let $q_i(i=1,\hdots,r)$ be any invertible 
elements of $R$ such that $q_i=q_j$ if $s_i$ is conjugate to $s_j$ in $W$. The Iwahori-Hecke algebra 
$\mathcal{H}_{R}(W,S)$ is an $R$-algebra generated by $\{T_{s_i}|s_i{\in}S\}$ with the relations: 
\renewcommand{\theenumi}{\arabic{enumi}}
\renewcommand{\labelenumi}{(H\theenumi)}
\begin{enumerate}
\item
$T_{s_i}^2 = (q_i-q_i^{-1})T_{s_i}+1$ \qquad if $i=1,2,\hdots,r$,
\item
$(T_{s_i}T_{s_j})^{k_{ij}}=(T_{s_j}T_{s_i})^{k_{ij}}$ \qquad if $m_{ij}=2k_{ij}$,
\item
$(T_{s_i}T_{s_j})^{k_{ij}}T_{s_i}=(T_{s_j}T_{s_i})^{k_{ij}}T_{s_j}$ \qquad if $m_{ij}=2k_{ij}+1$,
\end{enumerate}
\renewcommand{\theenumi}{\arabic{enumi}}
\renewcommand{\labelenumi}{(\theenumi)}
where $m_{ij}$ is the order of $s_is_j$ in $W$. 
We define $T_w=T_{s_{i_1}}T_{s_{i_2}}{\cdots}T_{s_{i_k}}$ where $w=s_{i_1}s_{i_2}{\cdots}s_{i_k}$ is 
a reduced expression of $w$. It is known that $T_w$ is well defined because 
two elements $T_w$ and $T_{w'}$, where $w$ and $w'$ are reduced expressions of 
an element of $W$, coincide and that $\{T_w|w{\in}W\}$ form a basis of $\mathcal{H}_{R}(W,S)$ as 
free $R$-modules. 
The relations (H1)--(H3) is equivalent to the following two relations: 
\renewcommand{\theenumi}{\arabic{enumi}}
\renewcommand{\labelenumi}{(h\theenumi)}
\begin{enumerate}
\item
$T_{s_i}T_w = T_{s_iw}$ \qquad if $l(w)<l(s_iw)$,
\item
$T_{s_i}T_w = (q_i-q_i^{-1})T_w + T_{s_iw}$ \qquad if $l(w)>l(s_iw)$,
\end{enumerate}
or equivalently, 
\renewcommand{\theenumi}{\arabic{enumi}}
\renewcommand{\labelenumi}{(h'\theenumi)}
\begin{enumerate}
\item
$T_{w}T_{s_i} = T_{ws_i}$ \qquad$ if l(w)<l(ws_i)$,
\item
$T_wT_{s_i} = (q_i-q_i^{-1})T_w + T_{ws_i}$ \qquad if $l(w)>l(ws_i)$,
\end{enumerate}
\renewcommand{\theenumi}{\arabic{enumi}}
\renewcommand{\labelenumi}{(\theenumi)}
where $l(w)$ means the length of $w$. 
We write $T_i=T_{s_i}$ for brevity. \par
If $(W,S)$ is of type $A$ and of rank $r-1$, then $W$ is isomorphic to the symmetric group 
$\mathfrak{S}_r$. 
Furthermore, all the elements of $S$ are conjugate to each other, hence we may assume 
$q_1=\cdots=q_{r-1}=q$. 
The Iwahori-Hecke algebra $\mathcal{H}_{R,r}(q)=\mathcal{H}_{R}(W,S)$ of type $A$ has 
defining relations: 
\renewcommand{\theenumi}{\arabic{enumi}}
\renewcommand{\labelenumi}{(A\theenumi)}
\begin{enumerate}
\item
$T_i^2 = (q-q^{-1})T_i+1$ \qquad if $i=1,2,\hdots,r-1$,
\item
$T_iT_{i+1}T_i = T_{i+1}T_iT_{i+1}$ \qquad if $i=1,2,\hdots,r-2$,
\item
$T_iT_j = T_jT_i$ \qquad if $|i-j|>1$.
\end{enumerate}
\renewcommand{\theenumi}{\arabic{enumi}}
\renewcommand{\labelenumi}{(\theenumi)}
Let $V={\oplus}_{k=1}^{m+n}Rv_k$ be a $\mathbb{Z}_2$-graded $R$-module of rank $m+n$. 
By $\mathbb{Z}_2$-graded, we mean that $V$ is a direct sum of two submodules 
$V_{\bar{0}}=\oplus_{k=1}^{m}Rv_k$ and $V_{\bar{1}}=\oplus_{k=m+1}^{m+n}Rv_k$, and that 
for each homogeneous element the degree map $|{\cdot}|$ 
\begin{equation*}
|v|=
\begin{cases}
0& \text{if $v{\in}V_{\bar{0}}$},\\
1& \text{if $v{\in}V_{\bar{1}}$,}
\end{cases}
\end{equation*}
is given. 
In order to define a representation of $\mathcal{H}_{R,r}(q)$ on $V^{{\otimes}r}$, 
we define a right operator $T$ on $V{\otimes}V$ as follows. 
\begin{equation}
v_k{\otimes}v_lT=
\begin{cases}
\dfrac{(-1)^{|v_k|}(q+q^{-1})+q-q^{-1}}{2}v_k{\otimes}v_l& \text{if $k=l$,}\\
(-1)^{|v_k||v_l|}v_l{\otimes}v_k+(q-q^{-1})v_k{\otimes}v_l& \text{if $k<l$,}\\
(-1)^{|v_k||v_l|}v_l{\otimes}v_k& \text{if $k>l$.}
\end{cases}
\end{equation}
The factor $2^{-1}$ in the case $k=l$ vanishes whether $v_k{\in}V_{\bar{0}}$ or 
$v_k{\in}V_{\bar{1}}$. So the above definition makes sense on $R$. 
Now we obtain right operators 
$\operatorname{Id}^{{\otimes}i-1}{\otimes}T{\otimes}\operatorname{Id}^{{\otimes}r-i-1}$
($i=1,2,\hdots,r-1$) 
on $V^{{\otimes}r}$ where $\operatorname{Id}$ is the identity operator on $V$. 
Let us define a map 
$\pi_r:\{T_1,\hdots,T_{r-1}\}{\longrightarrow}\operatorname{End}_R(V^{{\otimes}r})$ by 
$\pi_r(T_i)=\operatorname{Id}^{{\otimes}i-1}{\otimes}T{\otimes}\operatorname{Id}^{{\otimes}r-i-1}$
($i=1,2,\hdots,r-1$). 
\begin{Theo}
$\pi_r$ defines a representation of $\mathcal{H}_{R,r}(q)$ on $V^{{\otimes}r}$. 
\end{Theo}
\begin{proof}
One can check that the above operators satisfy the defining relations (A1)--(A3) 
by a direct computation. For example, the relation (A1) is shown as follows. \\
case1 : $k=l$
\begin{equation*}
\begin{split}
v_k{\otimes}v_kT^2&=\dfrac{1}{4}\big{\{}(-1)^{|v_k|}(q+q^{-1})+q-q^{-1}\big{\}}^2v_k{\otimes}v_k\\
&=\dfrac{1}{2}\big{\{}q^2+q^{-2}+(-1)^{|v_k|}(q+q^{-1})(q-q^{-1})\big{\}}v_k{\otimes}v_k
\end{split}
\end{equation*}
\begin{equation*}
\begin{split}
v_k{\otimes}v_k\big{\{}(q-q^{-1})T+1\big{\}}&=
\dfrac{(q-q^{-1})}{2}\big{\{}(-1)^{|v_k|}(q+q^{-1})+q-q^{-1}\big{\}}v_k{\otimes}v_k+v_k{\otimes}v_k\\
&=\dfrac{1}{2}\big{\{}q^2+q^{-2}+(-1)^{|v_k|}(q+q^{-1})(q-q^{-1})\big{\}}v_k{\otimes}v_k
\end{split}
\end{equation*}
case2 : $k<l$
\begin{equation*}
\begin{split}
v_k{\otimes}v_lT^2&=\big{\{}(-1)^{|v_k||v_l|}v_l{\otimes}v_k+(q-q^{-1})v_k{\otimes}v_l\big{\}}T\\
&=v_k{\otimes}v_l+(q-q^{-1})
\big{\{}(-1)^{|v_k||v_l|}v_l{\otimes}v_k+(q-q^{-1})v_k{\otimes}v_l\big{\}}\\
&=(q-q^{-1})(-1)^{|v_k||v_l|}v_l{\otimes}v_k+(q^2+q^{-2}-1)v_k{\otimes}v_l
\end{split}
\end{equation*}
\begin{equation*}
\begin{split}
v_k{\otimes}v_l\big{\{}(q-q^{-1})T+1\big{\}}&=
(q-q^{-1})\big{\{}(-1)^{|v_k||v_l|}v_l{\otimes}v_k+(q-q^{-1})v_k{\otimes}v_l\big{\}}
+v_k{\otimes}v_l\\
&=(q-q^{-1})(-1)^{|v_k||v_l|}v_l{\otimes}v_k+(q^2+q^{-2}-1)v_k{\otimes}v_l
\end{split}
\end{equation*}
case3 : $k>l$
\begin{equation*}
\begin{split}
v_k{\otimes}v_lT^2&=
(-1)^{|v_k||v_l|}\big{\{}(-1)^{|v_k||v_l|}v_k{\otimes}v_l+(q-q^{-1})v_l{\otimes}v_k\big{\}}\\
&=v_k{\otimes}v_l+(-1)^{|v_k||v_l|}(q-q^{-1})v_l{\otimes}v_k
\end{split}
\end{equation*}
\begin{equation*}
\begin{split}
v_k{\otimes}v_l\big{\{}(q-q^{-1})T+1\big{\}}&=
(q-q^{-1})(-1)^{|v_k||v_l|}v_l{\otimes}v_k+v_k{\otimes}v_l
\end{split}
\end{equation*}
(A2) can be shown in a similar manner to (A1), albeit slightly lengthy. 
The relation (A3) is obvious. 
\end{proof}
\begin{Rem}
Another definition of the Iwahori-Hecke algebra, which is frequently used 
\renewcommand{\theenumi}{\arabic{enumi}}
\renewcommand{\labelenumi}{(A'\theenumi)}
\begin{enumerate}
\item
${\T}_i^2 = (q'-1){\T}_i+q'$ \qquad if $i=1,2,\hdots,r-1$,
\item
${\T}_i{\T}_{i+1}{\T}_i = {\T}_{i+1}{\T}_i{\T}_{i+1}$ \qquad if $i=1,2,\hdots,r-2$,
\item
${\T}_i{\T}_j = {\T}_j{\T}_i$ \qquad if $|i-j|>1$,
\end{enumerate}
\renewcommand{\theenumi}{\arabic{enumi}}
\renewcommand{\labelenumi}{(\theenumi)}
can be obtained from previous definition (A1)--(A3) by letting $q'=q^2$ and ${\T}_i=qT_i$. 
Accordingly, we get another right operator $\T$ on $V{\otimes}V$ as follows. 
\begin{equation}
v_k{\otimes}v_l{\T}=
\begin{cases}
\dfrac{(-1)^{|v_k|}(q'+1)+q'-1}{2}v_k{\otimes}v_k& \text{if $k=l$,}\\
(-1)^{|v_k||v_l|}\sqrt{q'}v_l{\otimes}v_k+(q'-1)v_k{\otimes}v_l& \text{if $k<l$,}\\
(-1)^{|v_k||v_l|}\sqrt{q'}v_l{\otimes}v_k& \text{if $k>l$.}
\end{cases}
\end{equation}
We can also define the action of the alternative generators ${\T}_i$ of the Iwahori-Hecke algebra 
in a similar manner. 
\end{Rem}
This representation $\pi_r$ is reduced to the (normal) $q$-permutation representation of 
$\mathcal{H}_{R,r}(q)$ obtained from Drinfeld-Jimbo solutions to the Yang-Baxter equation 
when $n=0$ and to the sign permutation action (see \cite{B-R}) 
of the symmetric group when $q{\rightarrow}1$. 

\section{The quantum superalgebra $U_q^\sigma\big{(}{\gl}(m,n)\big{)}$ and the vector representation}
In Kac's paper\cite{Kac}, classical superalgebras have been classified and studied in detail. 
Quantum superalgebras have been defined in several articles such as \cite{B-K-K},\cite{K-T} 
or \cite{Yamane}. Each definition of them seems to be based on \cite{Kac} essentially. 
$U_q^\sigma\big{(}{\gl}(m,n)\big{)}$ is a Hopf algebra obtained from the \lq\lq naive" quantum 
superalgebra $U_q\big{(}{\gl}(m,n)\big{)}$, which is a Hopf superalgebra, 
by adding an involutive element $\sigma$. 
According to \cite{B-K-K}, we adopt $U_q^\sigma\big{(}{\gl}(m,n)\big{)}$ to construct 
the vector representation on the tensor space $V^{{\otimes}r}$. 
\par
Let $\Pi=\{\alpha_i\}_{i{\in}I}$ be a set of simple roots with the index set $I=\{1,\hdots,r\}$. 
We assume that $I$ is a disjoint union of two subsets $I_{\rm even}$ and $I_{\rm odd}$. 
We define a map $p:I{\longrightarrow}\{0,1\}$ to be such that
\begin{equation*}
p(i)=
\begin{cases}
0& \text{if $i{\in}I_{\rm even}$},\\
1& \text{if $i{\in}I_{\rm odd}$.}
\end{cases}
\end{equation*}
Let $P$ be a free $\mathbb{Z}$-module which includes all $\alpha_i{\in}P$($i{\in}I$). 
We assume that a $\mathbb{Q}$-valued symmetric bilinear form on $P$ 
$(\cdot,\cdot) : P{\times}P{\longrightarrow}{\mathbb{Q}}$ is defined and that the 
simple coroots $h_i{\in}P^*$($i{\in}I$) are given as data. The natural pairing 
${\langle}\cdot,\cdot{\rangle} : P^*{\times}P{\longrightarrow}{\mathbb{Z}}$ 
between $P$ and $P^*$ is assumed to satisfy 
\begin{equation*}
{\langle}h_i,\alpha_j{\rangle}=
\begin{cases}
2& \text{if $i=j$ and $i{\in}I_{\rm even}$,}\\
0\quad\text{or}\quad2& \text{if $i=j$ and $i{\in}I_{\rm odd}$,}\\
{\leq}0& \text{if $i{\neq}j$.}
\end{cases}
\end{equation*}
We denote by $\Pi^{\vee}=\{h_i|i{\in}I\}$ the set of all coroots. 
Furthermore, for each $i{\in}I$ we assume that there exists a nonzero integer $\ell_i$ such that 
$\ell_i{\langle}h_i,\lambda{\rangle}=(\alpha_i,\lambda)$ for every $\lambda{\in}P$.
Then we immediately have the Cartan matrix $A=[{\langle}h_i,\alpha_j{\rangle}]_{ij}$ 
is symmetrizable because $\ell_i{\langle}h_i,\alpha_j{\rangle}=(\alpha_i,\alpha_j)
=(\alpha_j,\alpha_i)=\ell_j{\langle}h_j,\alpha_i{\rangle}$. We mention that 
the symmetrized matrix is $A^{\rm sym}=\operatorname{diag}(\ell_1,\hdots,\ell_r)A=
[(\alpha_i,\alpha_j)]_{ij}$. Let $\mathfrak{h}=P^*{\otimes}_{\mathbb{Z}}\mathbb{Q}$. 
Then $\Phi=(\mathfrak{h},\Pi^{\vee},\Pi)$ 
is said to be a fundamental root data associated to $A$. 
Let $\g=\g(\Phi)$ be the contragredient Lie superalgebra obtained from $\Phi$ and $p$. 
According to \cite{B-K-K}, we define 
the quantized enveloping algebra $U_q(\g)$ to be the unital associative algebra over 
$\mathbb{Q}(q)$ with generators $q^h (h{\in}P^*),e_i,f_i (i{\in}I)$, which satisfy the 
following defining relations (compare \cite{K-T} and \cite{Yamane}):
\renewcommand{\theenumi}{\arabic{enumi}}
\renewcommand{\labelenumi}{(Q\theenumi)}
\begin{enumerate}
\item
$q^h=1$ \quad for $h=0$, 
\item
$q^{h_1}q^{h_2}=q^{h_1+h_2}$ \quad for $h_1,h_2{\in}P^*$, 
\item
$q^he_i=q^{{\langle}h,\alpha_j{\rangle}}e_iq^h$ \quad for $h{\in}P^*$ and $i{\in}I$,
\item
$q^hf_i=q^{-{\langle}h,\alpha_j{\rangle}}f_iq^h$ \quad for $h{\in}P^*$ and $i{\in}I$,
\item
$[e_i,f_j]=\delta_{ij}\dfrac{q^{\ell_ih_i}-q^{-\ell_ih_i}}{q^{\ell_i}-q^{-\ell_i}}$ 
\quad for $i,j{\in}I$,
\end{enumerate}
\renewcommand{\theenumi}{\arabic{enumi}}
\renewcommand{\labelenumi}{(\theenumi)}
where $[e_i,f_j]$ means the supercommutator 
\begin{equation*}
[e_i,f_j]=e_if_j-(-1)^{p(i)p(j)}f_je_i.
\end{equation*}
We assume further conditions (bitransitivity condition, see \cite{Kac} p.19):
\renewcommand{\theenumi}{\arabic{enumi}}
\renewcommand{\labelenumi}{(Q\theenumi)}
\begin{enumerate}
\item[(Q6)]
If $a{\in}\sum_{i{\in}I}U_q(\mathfrak{n}_+)e_iU_q(\mathfrak{n}_+)$ satisfies 
$f_ia{\in}U_q(\mathfrak{n}_+)f_i$ for all $i{\in}I$, then $a=0$, 
\item[(Q7)]
If $a{\in}\sum_{i{\in}I}U_q(\mathfrak{n}_-)f_iU_q(\mathfrak{n}_-)$ satisfies 
$e_ia{\in}U_q(\mathfrak{n}_-)e_i$ for all $i{\in}I$, then $a=0$, 
\end{enumerate}
\renewcommand{\theenumi}{\arabic{enumi}}
\renewcommand{\labelenumi}{(\theenumi)}
where $U_q(\mathfrak{n}_+)$ (resp.$U_q(\mathfrak{n}_-)$) is the subalgebra of $U_q(\g)$ generated 
by $\{e_i|i{\in}I\}$ (resp. $\{f_i|i{\in}I\}$). 
$U_q(\g)$ is a Hopf superalgebra whose comultiplication $\triangle$, counit $\varepsilon$, 
antipode $S$ are as follows. 
\begin{equation*}
\begin{split}
&\triangle(q^h)=q^h{\otimes}q^h \quad \text{for $h{\in}P^*$}, \\
&\triangle(e_i)=e_i{\otimes}q^{-\ell_ih_i}+1{\otimes}e_i \quad \text{for $i{\in}I$}, \\
&\triangle(f_i)=f_i{\otimes}1+q^{\ell_ih_i}{\otimes}f_i \quad \text{for $i{\in}I$}, \\
&\varepsilon(q^h)=1 \quad \text{for $h{\in}P^*$}, \quad 
\varepsilon(e_i)=\varepsilon(f_i)=0 \quad \text{for $i{\in}I$}, \\
&S(q^{{\pm}h})=q^{{\mp}h} \quad \text{for $h{\in}P^*$}, \\
&S(e_i)=-e_iq^{\ell_ih_i},{\quad}S(f_i)=-q^{-\ell_ih_i}f_i \quad \text{for $i{\in}I$}. 
\end{split}
\end{equation*}
This is not a Hopf algebra. 
In order to give a Hopf algebra structure to $U_q(\g)$, we define an involutive operator 
$\sigma$ on $U_q(\g)$ by $\sigma(q^h)=q^h$ for all $h{\in}P^*$ and $\sigma(e_i)=(-1)^{p(i)}e_i$,
$\sigma(f_i)=(-1)^{p(i)}f_i$ for all $i{\in}I$. 
Let $U^{\sigma}_q(\g)=U_q(\g){\oplus}U_q(\g)\sigma$. Then $U^{\sigma}_q(\g)$ is the algebra 
with the additional multiplication law given by $\sigma^2=1$ and $\sigma^{-1}x\sigma=\sigma(x)$ for 
any $x{\in}U_q(\g)$. 
$U_q(\g)$ is a Hopf algebra whose comultiplication $\triangle_{\sigma}$, 
counit $\varepsilon_{\sigma}$, antipode $S_{\sigma}$ are as follows. 
\begin{equation*}
\begin{split}
&\triangle_{\sigma}(\sigma)=\sigma{\otimes}\sigma, \\
&\triangle_{\sigma}(q^h)=q^h{\otimes}q^h \quad \text{for $h{\in}P^*$}, \\
&\triangle_{\sigma}(e_i)=e_i{\otimes}q^{-\ell_ih_i}+\sigma^{p(i)}{\otimes}e_i \quad \text{for $i{\in}I$}, \\
&\triangle_{\sigma}(f_i)=f_i{\otimes}1+\sigma^{p(i)}q^{\ell_ih_i}{\otimes}f_i \quad \text{for $i{\in}I$}, \\
&\varepsilon_{\sigma}(\sigma)=\varepsilon_{\sigma}(q^h)=1 
\quad \text{for $h{\in}P^*$}, \quad 
\varepsilon_{\sigma}(e_i)=\varepsilon_{\sigma}(f_i)=0 \quad \text{for $i{\in}I$}, \\
&S_{\sigma}(\sigma)=\sigma, \quad S_{\sigma}(q^{{\pm}h})=q^{{\mp}h} 
\quad \text{for $h{\in}P^*$}, \\
&S_{\sigma}(e_i)=-\sigma^{p(i)}e_iq^{\ell_ih_i},{\quad}S_{\sigma}(f_i)=-\sigma^{p(i)}q^{-\ell_ih_i}f_i 
\quad \text{for $i{\in}I$}. 
\end{split}
\end{equation*}
The quantized enveloping algebra $U_q\big{(}\gl(m,n)\big{)}$ 
is obtained from the fundamental root data as follows. 
\begin{itemize}
\item
$I=I_{\rm even}{\cup}I_{\rm odd}$ is defined by 
$I_{\rm even}=\{1,2,\hdots,m-1,m+1,\hdots,m+n-1\}$ and $I_{\rm odd}=\{m\}$,
\item
$P=\oplus_{b{\in}B}\mathbb{Z}\epsilon_b$, where 
$B=B_+{\cup}B_-$ with $B_+=\{1,\hdots,m\}$ and $B_-=\{m+1,\hdots,m+n\}$,
\item
$(\cdot,\cdot):P{\times}P{\longrightarrow}\mathbb{Q}$ is 
the symmetric bilinear form on $P$ defined by
\[
(\epsilon_a,\epsilon_{a'})=
\begin{cases}
1& \text{if $a=a'{\in}B_+$,}\\
-1& \text{if $a=a'{\in}B_-$,}\\
0& \text{otherwise,}
\end{cases}
\]
\item
$\Pi=\{\alpha_i|i{\in}I\}$ is defined by $\alpha_i=\epsilon_i-\epsilon_{i+1}$,
\item
$\Pi^{\vee}=\{h_i|i{\in}I\}$ is uniquely determined by the formula 
$\ell_i{\langle}h_i,\lambda{\rangle}=(\alpha_i,\lambda)$ for any $\lambda{\in}P$, \\
where 
\[
\ell_i=
\begin{cases}
1& \text{if $i=1,\hdots,m$,}\\
-1& \text{if $i=m+1,\hdots,m+n-1$.}
\end{cases}
\]
\end{itemize}Let $V$ be as in section2, and suppose $R=\mathbb{Q}(q)$. 
The vector representation ($\rho,V$) of $U^\sigma_q\big{(}\gl(m,n)\big{)}$ on $\mathbb{Z}_2$-graded 
vector space $V=V_{\bar{0}}{\oplus}V_{\bar{1}}$ (recall that 
$V_{\bar{0}}=\oplus_{i=1}^{m}Rv_i,V_{\bar{1}}=\oplus_{i=m+1}^{m+n}Rv_i$) is defined by 
\begin{equation}
\begin{split}
&\rho(\sigma)v_j=(-1)^{|v_j|}v_j \quad \text{for $j=1,\hdots,m+n$}, \\
&\rho(q^h)v_j=q^{\epsilon_j(h)}v_j \quad \text{for $h{\in}P^*,j=1,\hdots,m+n$}, \\
&\rho(e_j)v_{j+1}=v_j \quad \text{for $j=1,\hdots,m+n-1$}, \\
&\rho(f_j)v_j=v_{j+1} \quad \text{for $j=1,\hdots,m+n-1$}, \\
&\text{otherwise $0$.}
\end{split}
\end{equation}
This representation can be extended to the representation on the tensor space $V^{{\otimes}r}$. 
Let ${\rho}_r$ be the map from $U^\sigma_q\big{(}\gl(m,n)\big{)}$ to 
$\operatorname{End}_{\mathbb{Q}(q)}(V^{{\otimes}r})$ defined by 
\begin{equation}
\begin{split}
&{\rho}_r(\sigma)=\rho(\sigma)^{{\otimes}r}, \\
&{\rho}_r(q^h)=\rho(q^h)^{{\otimes}r} \quad \text{for $h{\in}P^*$}, \\
&{\rho}_r(e_i)=\sum_{k=1}^N\rho(\sigma^{p(i)})^{{\otimes}k-1}{\otimes}\rho(e_i) 
{\otimes}\rho(q^{-\ell_ih_i})^{{\otimes}r-k} \quad 
\text{for $i{\in}I$}, \\
&{\rho}_r(f_i)=\sum_{k=1}^r\rho(\sigma^{p(i)}q^{\ell_ih_i})^{{\otimes}k-1}{\otimes}\rho(f_i)
{\otimes}\operatorname{Id}^{{\otimes}r-k} \quad 
\text{for $i{\in}I$}.
\end{split}
\end{equation}
\begin{Prop}[\cite{B-K-K} Proposition3.1]
${\rho}_r$ gives a completely reducible representation of $U^\sigma_q\big{(}\gl(m,n)\big{)}$ 
on $V^{{\otimes}r}$ for $r{\ge}1$. 
\end{Prop}
Making use of the comultiplication of Hopf algebra, we obtain ${\rho}_r$. 
Let $\triangle^{(1)}=\triangle_\sigma$ at first and set 
$\triangle^{(k)}=(\triangle_\sigma{\otimes}\operatorname{Id}^{{\otimes}k-1})\triangle^{(k-1)}$ 
inductively. Then from the definition of $\triangle_\sigma$, we have 
${\rho}_r(x)=\rho^{{\otimes}r}\circ\triangle^{(r-1)}(x)$ for $x{\in}U^\sigma_q\big{(}\gl(m,n)\big{)}$ 
immediately. 

\section{Schur-Weyl reciprocity for $U^\sigma_q\big{(}\gl(m,n)\big{)}$ and 
$\mathcal{H}_{\mathbb{Q}(q),r}(q)$}
In the preceding sections, we obtained the left action of $U^\sigma_q\big{(}\gl(m,n)\big{)}$ 
and the right one of $\mathcal{H}_{\mathbb{Q}(q),r}(q)$ on the tensor space. 
We notice that the element $q$ of $\mathbb{Q}(q)$ is an indeterminate. 
In this section we derive the commutativity between these two actions. 
Furthermore, 
we establish Schur-Weyl reciprocity for these two algebras. 
We consider the relation between the vector representation ${\rho}_2$ and the operator $T$, 
both act on the tensor space $V^{{\otimes}2}$, at first. 
\begin{Prop}
$T$ commutes with the action of $U^\sigma_q\big{(}\gl(m,n)\big{)}$ on $V^{{\otimes}2}$ which is 
given by the representation $\rho_2$. 
\end{Prop}
\begin{proof}
The action $\rho_2(g)=(\rho{\otimes}\rho)\circ\triangle_\sigma(g)$ of 
$U^\sigma_q\big{(}\gl(m,n)\big{)}$ on $V^{{\otimes}2}$ is defined in 
several cases depending upon the generator of $U^\sigma_q\big{(}\gl(m,n)\big{)}$ and the basis 
vector $v_i{\otimes}v_j$ of $V^{{\otimes}2}$. \\
Case1 : $g=\sigma$\\
One can immediately check the commutativity between $\rho_2$ and $T$. \\
Case2 : $g=q^h$($h{\in}P^*$)\\
It is clear in this case because of the equation, 
\[
\rho_2(q^h)v_i{\otimes}v_j
=q^{\epsilon_i(h)}q^{\epsilon_j(h)}v_i{\otimes}v_j
=q^{(\epsilon_i+\epsilon_j)(h)}v_i{\otimes}v_j.
\]
Case3 : $g=e_k$ ($k=1,2,\hdots,m+n$) \\
We obtain the case-by-case definition of $\rho_2$ from (3.2) as follows. 
\begin{equation}
\rho_2(e_k)v_i{\otimes}v_j=
\begin{cases}
0& \text{if $i,j{\neq}k+1$},\\
v_{i-1}{\otimes}v_j& \text{if $i=k+1,j{\neq}k,k+1$},\\
q^{-(-1)^{|v_j|}}v_{i-1}{\otimes}v_j& \text{if $i=k+1,j=k$},\\
q^{(-1)^{|v_j|}}v_{i-1}{\otimes}v_j\\
 +(-1)^{|v_i|\delta_{km}}v_i{\otimes}v_{j-1}& 
\text{if $i=k+1,j=k+1$},\\
(-1)^{|v_i|\delta_{km}}v_i{\otimes}v_{j-1}& \text{if $i{\neq}k+1,j=k+1$}.\\
\end{cases}
\end{equation}
The operator $T$ has already defined in (2.1). 
One can check the commutativity between $T$ and $\rho_2$ by a direct computation in each case. \\
Case3--1 : $i,j{\neq}k+1$ \\
In this case, we imediately see 
$\big{(}\rho_2(e_k)v_i{\otimes}v_j\big{)}T
=\rho_2(e_k)(v_i{\otimes}v_jT)=0$. \\
Case3--2 : $i=k+1,j>k+1$ 
\begin{equation*}
\begin{split}
\big{\{}\rho_2(e_k)v_i{\otimes}v_j\big{\}}T
&=v_{i-1}{\otimes}v_jT \\
&=(-1)^{|v_{i-1}||v_j|}v_j{\otimes}v_{i-1}+(q-q^{-1})v_{i-1}{\otimes}v_j \\
\rho_2(e_k)(v_i{\otimes}v_jT)
&=\rho_2(e_k)\big{\{}(-1)^{|v_i||v_j|}v_j{\otimes}v_i+(q-q^{-1})v_i{\otimes}v_j\big{\}} \\
&=(-1)^{(|v_i|+\delta_{km})|v_j|}v_j{\otimes}v_{i-1}+(q-q^{-1})v_{i-1}{\otimes}v_j
\end{split}
\end{equation*}
If $k<m$, then $|v_{i-1}|=|v_i|=\delta_{km}=0$. 
If $k=m$, then $|v_{i-1}|=0,|v_i|=1,\delta_{km}=1$. 
If $k>m$, then $|v_{i-1}|=|v_i|=1,\delta_{km}=0$. 
In each case, $\big{\{}\rho_2(e_k)v_i{\otimes}v_j\big{\}}T=\rho_2(e_k)(v_i{\otimes}v_jT)$ holds. \\
Case3--3 : $i=k+1,j<k$ 
\begin{equation*}
\begin{split}
\big{\{}\rho_2(e_k)v_i{\otimes}v_j\big{\}}T
&=v_{i-1}{\otimes}v_jT \\
&=(-1)^{|v_{i-1}||v_j|}v_j{\otimes}v_{i-1} \\
\rho_2(e_k)(v_i{\otimes}v_jT)
&=\rho_2(e_k)\big{\{}(-1)^{|v_i||v_j|}v_j{\otimes}v_i\big{\}} \\
&=(-1)^{(|v_i|+\delta_{km})|v_j|}v_j{\otimes}v_{i-1}
\end{split}
\end{equation*}
As in the case3--2, $\big{\{}\rho_2(e_k)v_i{\otimes}v_j\big{\}}T=\rho_2(e_k)(v_i{\otimes}v_jT)$ 
holds. \\
Case3--4 : $i=k+1,j=k$ 
\begin{equation*}
\begin{split}
\big{\{}\rho_2(e_k)v_i{\otimes}v_j\big{\}}T
&=q^{-(-1)^{|v_j|}}v_{i-1}{\otimes}v_jT \\
&=q^{-(-1)^{|v_j|}}2^{-1}\big{\{}(-1)^{|v_j|}(q+q^{-1})+q-q^{-1}\big{\}}v_j{\otimes}v_j \\
\rho_2(e_k)(v_i{\otimes}v_jT)
&=\rho_2(e_k)\big{\{}(-1)^{|v_i||v_j|}v_j{\otimes}v_i\big{\}} \\
&=(-1)^{(|v_{j+1}|+\delta_{km})|v_j|}v_j{\otimes}v_j
\end{split}
\end{equation*}
If $k<m$, then $|v_j|=\delta_{km}=0$. 
If $k=m$, then $|v_j|=0,\delta_{km}=1$. 
If $k>m$, then $|v_j|=|v_{j+1}|=1,\delta_{km}=0$. 
In each case, $\big{\{}\rho_2(e_k)v_i{\otimes}v_j\big{\}}T=\rho_2(e_k)(v_i{\otimes}v_jT)$ holds. \\
Case3--5 : $i=k+1,j=k+1$ 
\begin{equation*}
\begin{split}
\big{\{}\rho_2(e_k)v_i{\otimes}v_j\big{\}}T
&=\big{\{}q^{(-1)^{|v_j|}}v_{i-1}{\otimes}v_j+(-1)^{|v_i|\delta_{km}}v_i{\otimes}v_{j-1}\big{\}}T \\
&=q^{(-1)^{|v_i|}}(-1)^{|v_{i-1}||v_i|}v_i{\otimes}v_{i-1}
+\big{\{}q^{(-1)^{|v_i|}}(q-q^{-1})+(-1)^{(|v_{i-1}|+\delta_{km})|v_i|}\big{\}}v_{i-1}{\otimes}v_i \\
\rho_2(e_k)(v_i{\otimes}v_jT)
&=\rho_2(e_k)\big{\{}\dfrac{(-1)^{|v_i|}(q+q^{-1})+q-q^{-1}}{2}v_i{\otimes}v_i\big{\}} \\
&=\big{\{}\dfrac{(-1)^{|v_i|}(q+q^{-1})+q-q^{-1}}{2}\big{\}}\big{\{}q^{(-1)^{|v_i|}}
v_{i-1}{\otimes}v_i+(-1)^{|v_i|\delta_{km}}v_i{\otimes}v_{i-1}\big{\}}
\end{split}
\end{equation*}
If $k<m$, then $|v_{i-1}|=|v_i|=\delta_{km}=0$. 
If $k=m$, then $|v_{i-1}|=0,|v_i|=\delta_{km}=1$. 
If $k>m$, then $|v_{i-1}|=|v_i|=1,\delta_{km}=0$. 
In each case, $\big{\{}\rho_2(e_k)v_i{\otimes}v_j\big{\}}T=\rho_2(e_k)(v_i{\otimes}v_jT)$ holds. \\
Case3--6 : $i<k,j=k+1$ 
\begin{equation*}
\begin{split}
\big{\{}\rho_2(e_k)v_i{\otimes}v_j\big{\}}T
&=(-1)^{|v_i|\delta_{km}}v_i{\otimes}v_{j-1}T \\
&=(-1)^{(|v_{j-1}|+\delta_{km})|v_i|}v_{j-1}{\otimes}v_i
+(-1)^{|v_i|\delta_{km}}(q-q^{-1})v_i{\otimes}v_{j-1} \\
\rho_2(e_k)(v_i{\otimes}v_jT)
&=\rho_2(e_k)\big{\{}(-1)^{|v_i||v_j|}v_j{\otimes}v_i+(q-q^{-1})v_i{\otimes}v_j\big{\}} \\
&=(-1)^{|v_i||v_j|}v_{j-1}{\otimes}v_i+(q-q^{-1})(-1)^{|v_i|\delta_{km}}v_i{\otimes}v_{j-1}
\end{split}
\end{equation*}
If $k<m$, then $|v_i|=|v_{j-1}|=|v_j|=\delta_{km}=0$. 
If $k=m$, then $|v_i|=|v_{j-1}|=0,|v_j|=\delta_{km}=1$. 
If $k>m$, then $|v_{j-1}|=|v_j|=1,\delta_{km}=0$. 
In each case, $\big{\{}\rho_2(e_k)v_i{\otimes}v_j\big{\}}T=\rho_2(e_k)(v_i{\otimes}v_jT)$ holds. \\
Case3-7 : $i=k,j=k+1$ 
\begin{equation*}
\begin{split}
\big{\{}\rho_2(e_k)v_i{\otimes}v_j\big{\}}T
&=(-1)^{|v_i|\delta_{km}}v_i{\otimes}v_{j-1}T \\
&=(-1)^{|v_i|\delta_{km}}\dfrac{(-1)^{|v_i|}(q+q^{-1})+q-q^{-1}}{2}v_i{\otimes}v_i \\
\rho_2(e_k)(v_i{\otimes}v_jT)
&=\rho_2(e_k)\big{\{}(-1)^{|v_i||v_j|}v_j{\otimes}v_i+(q-q^{-1})v_i{\otimes}v_j\big{\}} \\
&=\big{\{}(-1)^{|v_i||v_j|}q^{-(-1)|v_i|}+(q-q^{-1})(-1)^{|v_i|\delta_{km}}v_i{\otimes}v_i
\end{split}
\end{equation*}
If $k<m$, then $|v_i|=|v_j|=\delta_{km}=0$. 
If $k=m$, then $|v_i|=0,|v_j|=\delta_{km}=1$. 
If $k>m$, then $|v_i|=|v_j|=1,\delta_{km}=0$. 
In each case, $\big{\{}\rho_2(e_k)v_i{\otimes}v_j\big{\}}T=\rho_2(e_k)(v_i{\otimes}v_jT)$ holds. \\
Case3-8 : $i>k+1,j=k+1$ 
\begin{equation*}
\begin{split}
\big{\{}\rho_2(e_k)v_i{\otimes}v_j\big{\}}T
&=(-1)^{|v_i|\delta_{km}}v_i{\otimes}v_{j-1}T \\
&=(-1)^{(|v_{j-1}|+\delta_{km})|v_i|}v_{j-1}{\otimes}v_i \\
\rho_2(e_k)(v_i{\otimes}v_jT)
&=\rho_2(e_k)(-1)^{|v_i||v_j|}v_j{\otimes}v_i \\
&=(-1)^{|v_i||v_j|}v_{j-1}{\otimes}v_i
\end{split}
\end{equation*}
If $k<m$, then $|v_{j-1}|=|v_j|=\delta_{km}=0$. 
If $k=m$, then $|v_{j-1}|=0,|v_j|=\delta_{km}=1$. 
If $k>m$, then $|v_{j-1}|=|v_j|=1,\delta_{km}=0$. 
In each case, $\big{\{}\rho_2(e_k)v_i{\otimes}v_j\big{\}}T=\rho_2(e_k)(v_i{\otimes}v_jT)$ holds. \\
Case3-1 to Case3-8 exhaust the possible cases in (4.1), thus we have checked the commutativity for case3. \\
Case4 : $g=f_k$ ($k=1,2,\hdots,m+n$) \\
In a similar manner to case3, we obtain 
\begin{equation*}
\rho_2(f_k)v_i{\otimes}v_j=
\begin{cases}
0& \text{if $i,j{\neq}k$},\\
(-1)^{|v_i|\delta_{km}}v_i{\otimes}v_{j+1}& \text{if $i{\neq}k,k+1,j=k$},\\
v_{i+1}{\otimes}v_j\\
+(-1)^{|v_i|\delta_{km}}q^{(-1)^{|v_i|}}v_i{\otimes}v_{j+1}& 
\text{if $i=k,j=k$},\\
(-1)^{|v_i|\delta_{km}}q^{-(-1)^{|v_i|}}v_i{\otimes}v_{j+1}& \text{if $i=k+1,j=k$},\\
v_{i+1}{\otimes}v_j& \text{if $i=k,j{\neq}k$},\\
\end{cases}
\end{equation*}
and one can check for this case, so we omit the detail. \par
Finally, these exhaust the entirely possible cases, thus we have completed the proof. 
\end{proof}
\begin{Prop}
For every $g{\in}\mathcal{H}_{\mathbb{Q}(q),r}(q)$ and $x{\in}U^\sigma_q\big{(}\gl(m,n)\big{)}$, 
we have $\pi_r(g)\rho_r(x)=\rho_r(x)\pi_r(g)$. 
\end{Prop}
\begin{proof}
When $r=2$, we have already shown in proposition 4.1. 
It is enough to prove for $g{\in}\{T_1,\hdots,T_{r-1}\}$. 
We deduce from cocomutativity of $\triangle_\sigma$ that 
\begin{equation*}
\begin{split}
\triangle^{(r)}&=(\triangle_\sigma{\otimes}\operatorname{id}^{{\otimes}r-1})\triangle^{(r-1)} \\
&=(\operatorname{id}{\otimes}\triangle_\sigma{\otimes}\operatorname{id}^{{\otimes}r-2})\triangle^{(r-1)} \\
&=(\operatorname{id}^{{\otimes}2}{\otimes}\triangle_\sigma{\otimes}\operatorname{id}^{{\otimes}r-3})
\triangle^{(r-1)} \\
&=\hdots \\
&=(\operatorname{id}^{{\otimes}r-1}{\otimes}\triangle_\sigma)\triangle^{(r-1)}
\end{split}
\end{equation*}
where $\operatorname{id}$ is the identity operator on $U^\sigma_q\big{(}\gl(m,n)\big{)}$. 
For any $r$ with $r>2$, we may write 
$\triangle^{(r-2)}(x)=x_1{\otimes}{\cdots}{\otimes}x_{r-1}$ for some 
$x_1,\hdots,x_{r-1}{\in}U^\sigma_q\big{(}\gl(m,n)\big{)}$. 
Then applying the case $r=2$, we have the following. 
\begin{equation*}
\begin{split}
\pi_r(T_i)\rho_r(x)&=
\big{\{}\operatorname{Id}^{{\otimes}i-1}{\otimes}\pi_2(T_i)
{\otimes}\operatorname{Id}^{{\otimes}r-i-1}\big{\}}
\big{\{}\rho^{{\otimes}r}\circ\triangle^{(r-1)}(x)\big{\}} \\
&=\big{\{}\operatorname{Id}^{{\otimes}i-1}{\otimes}\pi_2(T_i){\otimes}
\operatorname{Id}^{{\otimes}r-i-1}\big{\}}
\big{\{}\rho^{{\otimes}r}\circ
(\operatorname{id}^{{\otimes}i-1}{\otimes}\triangle_\sigma{\otimes}
\operatorname{id}^{{\otimes}r-i-1})\triangle^{(r-2)}(x)\big{\}} \\
&=\big{\{}\otimes_{k=1}^{i-1}\rho(x_k)\big{\}}{\otimes}\pi_2(T_i)\triangle_\sigma(x_i){\otimes}
\big{\{}\otimes_{l=i+1}^{r-1}\rho(x_l)\big{\}} \\
&=\big{\{}\otimes_{k=1}^{i-1}\rho(x_k)\big{\}}{\otimes}\triangle_\sigma(x_i)\pi_2(T_i){\otimes}
\big{\{}\otimes_{l=i+1}^{r-1}\rho(x_l)\big{\}} \\
&=\big{\{}\rho^{{\otimes}r}\circ(\operatorname{id}^{{\otimes}i-1}{\otimes}\triangle_\sigma{\otimes}
\operatorname{id}^{{\otimes}r-i-1})\triangle^{(r-2)}(x)\big{\}}
\big{\{}\operatorname{Id}^{{\otimes}i-1}{\otimes}\pi_2(T_i){\otimes}
\operatorname{Id}^{{\otimes}r-i-1}\big{\}} \\
&=\rho_r(x)\pi_r(T_i) 
\end{split}
\end{equation*}
\end{proof}
We may define $\mathcal{G}_q$ to be the subalgebra of 
$\operatorname{End}_{R_0}\big{(}({R_0}^{m+n})^{{\otimes}r}\big{)}{\cong}
\operatorname{Mat}\big{(}(m+n)^r,{R_0}\big{)}$
generated by the set 
$\{\rho_r(\sigma),\rho_r(q^h),\rho_r(e_i)\rho_r(f_i)|h{\in}P^*,i{\in}I\}$ 
because all the matrix elements of those generators are in $R_0$ from (3.1) and (3.2). 
For the same reason we may also define $\mathcal{S}_q$ the one generated by 
$\{\pi_r(T_j)|j=1,\hdots,r-1\}$. 
Let us define two subalgebras of $\operatorname{Mat}\big{(}(m+n)^r,{R_0}\big{)}$ as follows. 
\begin{equation*}
\begin{split}
\Tilde{\mathcal{S}}_q&=\{X{\in}\operatorname{Mat}\big{(}(m+n)^r,{R_0}\big{)}|XY=YX \text{ for all } 
Y{\in}\mathcal{S}_q\} \\
\Tilde{\mathcal{G}}_q&=\{X{\in}\operatorname{Mat}\big{(}(m+n)^r,{R_0}\big{)}|XY=YX \text{ for all } 
Y{\in}\mathcal{G}_q\}
\end{split}
\end{equation*}
Let $R_0=\mathbb{Q}[q,q^{-1}]$ be the Laurent polynomial ring. 
Let us define the specialization to a nonzero complex number $t$ to be a ring homomorphism 
$\varphi_t : R_0{\longrightarrow}\mathbb{C}$ with the condition $\varphi_t(q)=t$. 
$\mathbb{C}$ becomes $(\mathbb{C},R_0)$-bimodule, with $R_0$ acting from the right via 
$\varphi_t$. 
Applying the spacialization $\varphi_t$, we obtain the specialized 
algebras $\mathcal{G}_t=\mathbb{C}{\otimes_{R_0}}\mathcal{G}_q$ and 
$\mathcal{S}_t=\mathbb{C}{\otimes_{R_0}}\mathcal{S}_q$ which are subalgebras of 
$\operatorname{Mat}\big{(}(m+n)^r,\mathbb{C}\big{)}$. 
We also have $\Tilde{\mathcal{G}}_t=\mathbb{C}{\otimes_{R_0}}\Tilde{\mathcal{G}}_q$ and 
$\Tilde{\mathcal{S}}_t=\mathbb{C}{\otimes_{R_0}}\Tilde{\mathcal{S}}_q$ immediately. 
\begin{Prop}
$\mathcal{G}_q=\Tilde{\mathcal{S}}_q$ and $\mathcal{S}_q=\Tilde{\mathcal{G}}_q$ hold. 
\end{Prop}
\begin{proof}
Since ${R_0}$ is a principal ideal domain, the submodules $\mathcal{G}_q$ and $\mathcal{S}_q$ 
of the free ${R_0}$-module $\operatorname{Mat}\big{(}(m+n)^r,{R_0}\big{)}$ are also free. 
Let $X_i(q){\in}\operatorname{Mat}\big{(}(m+n)^r,{R_0}\big{)}$($i=1,\hdots,N$) be a basis of 
$\mathcal{G}_q$ and $x_i^{k,l}(q){\in}{R_0}$ the $(k,l)$-entry of $X_i(q)$. 
Then we immediately have that the specialized elements 
$X_i(t)$($i=1,\hdots,N$) generate $\mathcal{G}_t$ and 
$\dim_\mathbb{C}\mathcal{G}_t{\leq}\operatorname{rank}_{R_0}\mathcal{G}_q$. 
Because $X_i(q)$ are linearly independent, 
$\sum_{i=1}^N\alpha_i(q)x_i^{k,l}(q)=0$ for $\alpha_1(q),\hdots,\alpha_N(q){\in}{R_0}$ and 
for all $k,l$ imply $\alpha_1(q)=\hdots=\alpha_N(q)=0$. 
We consider the system of linear equations as follows. 
\begin{equation*}
\begin{bmatrix}
x_1^{1,1}(q) & x_2^{1,1}(q) & \cdots & x_N^{1,1}(q) \\
x_1^{1,2}(q) & x_2^{1,2}(q) & \cdots & x_N^{1,2}(q) \\
&\cdots\cdots\cdots \\
x_1^{(m+n)^r,(m+n)^r-1}(q) & x_2^{(m+n)^r,(m+n)^r-1}(q) & \cdots & x_N^{(m+n)^r,(m+n)^r-1}(q) \\
x_1^{(m+n)^r,(m+n)^r}(q) & x_2^{(m+n)^r,(m+n)^r}(q) & \cdots & x_N^{(m+n)^r,(m+n)^r}(q)
\end{bmatrix}
\begin{bmatrix}
\alpha_1(q) \\
\alpha_2(q) \\
\cdots \\
\alpha_N(q)
\end{bmatrix}
=
\begin{bmatrix}
0 \\
0 \\
\cdots \\
0 \\
0
\end{bmatrix}
\end{equation*}
This has only the trivial solution. 
But applying the specialization $\varphi_t$ to the above system, we possibly obtain a 
non-trivial solution. 
If there exists a non-trivial solution, then $t$ must be a zero of some 
Laurent polynomial whose coefficients are in $\mathbb{Q}$. 
Therefore if $t$ is a transcendental number, the above system has only the trivial 
solution and hence 
$\dim_\mathbb{C}\mathcal{G}_t=\operatorname{rank}_{R_0}\mathcal{G}_q$. 
In the same manner, we also have that if $t$ is a transcendental number, 
then $\dim_\mathbb{C}\mathcal{S}_t=\operatorname{rank}_{R_0}\mathcal{S}_q$. \par
We shall show that $\operatorname{rank}_{R_0}\Tilde{\mathcal{S}}_q=
\dim_\mathbb{C}\Tilde{\mathcal{S}}_t$. One can readily see that 
$\operatorname{rank}_{R_0}\Tilde{\mathcal{S}}_q{\geq}\dim_\mathbb{C}\Tilde{\mathcal{S}}_t$ 
where $\Tilde{\mathcal{S}}_t$ is the specialized algebra of $\Tilde{\mathcal{S}}_q$. 
Assume that $X(q)=\big{(}x^{k,l}(q)\big{)}{\in}\operatorname{Mat}\big{(}(m+n)^r,{R_0}\big{)}
{\in}\Tilde{\mathcal{S}}_q$. 
Then $X(q)$ commutes with $\pi_r(T_i)$ for all $i$, hence the matrix elements $x^{k,l}(q)$
($k,l=1,\hdots,(m+n)^{{\otimes}r}$) satisfy linear equations of coefficients in ${R_0}$. 
In the same manner as $\mathcal{G}_q$, one can find that 
the commutativity condition 
turns out to be the condition of solubilities of certain Laurent polynomials of coefficients in 
$\mathbb{Q}$. Thus if $t$ is a transcendental number, only the trivial equation exists, so 
$\operatorname{rank}_{R_0}\Tilde{\mathcal{S}}_q=
\dim_\mathbb{C}\Tilde{\mathcal{S}}_t$ holds. 
In a similar manner, we also have that if $t$ is a transcendental number, then 
$\operatorname{rank}_{R_0}\Tilde{\mathcal{G}}_q=
\dim_\mathbb{C}\Tilde{\mathcal{G}}_t$. 
Now we fix a transcendental number $t$. 
In the preceding work such like \cite{B-R}, \cite{Ser}, it has already been shown that 
$\Tilde{\mathcal{S}}_1=\mathcal{G}_1$. 
From this fact and the inequality 
\[
\dim_{\mathbb{C}}\mathcal{G}_1{\leq}\dim_{\mathbb{C}}\mathcal{G}_t
{\leq}\dim_{\mathbb{C}}\Tilde{\mathcal{S}}_t{\leq}\dim_{\mathbb{C}}\Tilde{\mathcal{S}}_1,
\]
we deduce that $\dim_{\mathbb{C}}\mathcal{G}_t=\dim_{\mathbb{C}}\Tilde{\mathcal{S}}_t$. 
Because $\mathcal{G}_t{\subseteq}\Tilde{\mathcal{S}}_t$, we obtain 
$\mathcal{G}_t=\Tilde{\mathcal{S}}_t$. 
It is known that the specialized algebra 
$\mathcal{H}_{\mathbb{C},r}(t)=\mathbb{C}{\otimes}_{R_0}\mathcal{H}_{R_0,r}(q)$ 
is (split)semisimple. 
Therefore applying the double centralizer theorem, we also obtain 
$\mathcal{S}_t=\Tilde{\mathcal{G}}_t$. 
Using the properties, 
\begin{equation*}
\begin{split}
\operatorname{rank}_{R_0}\mathcal{G}_q&=\dim_\mathbb{C}\mathcal{G}_t, \\
\operatorname{rank}_{R_0}\mathcal{S}_q&=\dim_\mathbb{C}\mathcal{S}_t, \\
\operatorname{rank}_{R_0}\Tilde{\mathcal{G}}_q&=\dim_\mathbb{C}\Tilde{\mathcal{G}}_t, \\ 
\operatorname{rank}_{R_0}\Tilde{\mathcal{S}}_q&=\dim_\mathbb{C}\Tilde{\mathcal{S}}_t, 
\end{split}
\end{equation*}
which are already shown in the previous discussion, 
we readily see that $\mathcal{G}_q=\Tilde{\mathcal{S}}_q$ and $\mathcal{S}_q=\Tilde{\mathcal{G}}_q$. 
\end{proof}
We denote $\pi_r\big{(}\mathcal{H}_{\mathbb{Q}(q),r}(q)\big{)}$ by $\mathcal{A}_q$ and 
$\rho_r\big{(}U^\sigma_q\big{(}\gl(m,n)\big{)}\big{)}$ by $\mathcal{B}_q$. 
Then we have the following. 
\begin{Theo}
$\operatorname{End}_{\mathcal{B}_q}V^{{\otimes}r}
=\mathcal{A}_q$ and $\operatorname{End}_{\mathcal{A}_q}V^{{\otimes}r}
=\mathcal{B}_q$ hold. 
\end{Theo}
\begin{proof}
Obviously $\mathcal{A}_q\cong\mathcal{S}_q{\otimes}_{R_0}\mathbb{Q}(q)$ and 
$\mathcal{B}_q\cong\mathcal{G}_q{\otimes}_{R_0}\mathbb{Q}(q)$ as $\mathbb{Q}(q)$-algebras. 
From proposition4.3 we obtain that 
$\operatorname{End}_{\mathcal{B}_q}V^{{\otimes}r}
=\mathcal{A}_q$ and $\operatorname{End}_{\mathcal{A}_q}V^{{\otimes}r}
=\mathcal{B}_q$, and we have completed the proof. 
\end{proof}

\section{Decomposition of the tensor space}
Let $\overline{\mathbb{Q}(q)}$ be the algebraic closure of $\mathbb{Q}(q)$. 
We define 
$\bar{U}^\sigma_q\big{(}\gl(m,n)\big{)}=U^\sigma_q\big{(}\gl(m,n)\big{)}
{\otimes_{\mathbb{Q}(q)}}\overline{\mathbb{Q}(q)}$, 
$\bar{\mathcal{A}}_q=\mathcal{A}_q{\otimes_{\mathbb{Q}(q)}}\overline{\mathbb{Q}(q)}$, 
$\bar{\mathcal{B}}_q=\mathcal{B}_q{\otimes_{\mathbb{Q}(q)}}\overline{\mathbb{Q}(q)}$. 
Then, 
$\pi_r\big{(}\mathcal{H}_{\overline{\mathbb{Q}(q)},r}(q)\big{)}=\bar{\mathcal{A}}_q$ and 
$\rho_r\big{(}\bar{U}^\sigma_q\big{(}\gl(m,n)\big{)}\big{)}=\bar{\mathcal{B}}_q$ 
as $\overline{\mathbb{Q}(q)}$-algebras of operators on 
$\bar{V}^{{\otimes}r}=(V{\otimes_{\mathbb{Q}(q)}}\overline{\mathbb{Q}(q)})^{{\otimes}r}$. 
Theorem4.4 still holds when we exchange the base field from $\mathbb{Q}(q)$ to 
$\overline{\mathbb{Q}(q)}$, namely, 
$\operatorname{End}_{\bar{\mathcal{B}}_q}\bar{V}^{{\otimes}r}
=\bar{\mathcal{A}}_q$ and $\operatorname{End}_{\bar{\mathcal{A}}_q}\bar{V}^{{\otimes}r}
=\bar{\mathcal{B}}_q$. \par
We denote by $\Par(r)$ the set of all partitions of $r$. 
By the double centralizer theorem, there is a subset $\Gamma$ of $\Par(r)$ such that 
$\bar{\mathcal{A}}_q=\oplus_{\lambda{\in}\Gamma}\bar{\mathcal{A}}_{q,\lambda}$ where 
$\bar{\mathcal{A}}_{q,\lambda}$ is the Wedderburn component corresponding to 
the irreducible representation of $\mathcal{H}_{\overline{\mathbb{Q}(q)},r}(q)$ 
indexed by $\lambda$. Moreover, we also obtain the decomposition of 
$\mathcal{H}_{\overline{\mathbb{Q}(q)},r}(q){\otimes}
\bar{U}^\sigma_q\big{(}\gl(m,n)\big{)}$-modules, 
\begin{equation}
\bar{V}^{{\otimes}r}=\bigoplus_{\lambda{\in}\Gamma}H_{\lambda}{\otimes}V_{\lambda}, 
\end{equation}
where $H_{\lambda}$ is the irreducible representation of 
$\mathcal{H}_{\overline{\mathbb{Q}(q)},r}(q)$ indexed by $\lambda$, and 
$V_{\lambda}$ is the one of $\bar{U}^\sigma_q\big{(}\gl(m,n)\big{)}$ such that 
$V_{\lambda}{\ncong}V_{\mu}$ if $\lambda\neq\mu$. 
Our subject in this chapter is to determine $\Gamma$. \par
Let $H(m,n;r)=\{\lambda=(\lambda_1,\lambda_2,\hdots){\in}\Par(r)|\lambda_j{\leq}n 
\text{ if } j>m\}$. Diagrams of elements of $H(m,n;r)$ are exactly those contained 
by the $(m,n)$-hooks. 
We shall show that $\Gamma=H(m,n;r)$. 
\begin{Theo}
$\bar{\mathcal{A}}_q=\bigoplus_{\lambda{\in}H(m,n;r)}\bar{\mathcal{A}}_{q,\lambda}$. Hence 
$\bar{V}^{{\otimes}r}=\bigoplus_{\lambda{\in}H(m,n;r)}H_{\lambda}{\otimes}V_{\lambda}$ holds. 
\end{Theo}
\begin{proof}
When $q=1$, then Berele and Regev have already shown that 
\begin{Theo}[\cite{B-R}3.20 The Hook Theorem]
Let $F$ be an algebraic closed field of characteristic $0$ and $\rho$ 
the sign permutation representation on $V^{{\otimes}n}$ where $V$ is a $(k,l)$-dimensional 
vector space over $F$. Then 
\[
\rho(F[\mathfrak{S}_n])=\bigoplus_{\lambda{\in}H(k,l;n)}A_{\lambda}
{\cong}\bigoplus_{\lambda{\in}H(k,l;n)}I_{\lambda}
\]
where each $A_{\lambda}$ is the Wedderburn component corresponding to $\lambda$, 
and $I_{\lambda}$ is a simple subalgebra of $F[\mathfrak{S}_n]$ such that 
$\rho(I_{\lambda})=A_{\lambda}$. 
\end{Theo}
Thus $\Gamma=H(m,n;r)$ holds for $q=1$. 
Let $t$ be a transcendental number. 
We have already shown in the proof of proposition4.3 that 
$\dim_{\mathbb{C}}\mathcal{G}_t=\dim_{\mathbb{C}}\mathcal{G}_1$ and 
$\dim_{\mathbb{C}}\mathcal{S}_t=\dim_{\mathbb{C}}\mathcal{S}_1$. 
Let $\mathcal{S}_t=\bigoplus_{\lambda{\in}\Par(r)}\mathcal{S}_{t,\lambda}$ 
be the Wedderburn decomposition. 
Then, by Theorem5.2, we have $\mathcal{S}_{1,\lambda}=0$ if and only if 
$\lambda{\not\in}H(m,n;r)$. 
Because $\dim_{\mathbb{C}}\mathcal{S}_{t,\lambda}{\geq}
\dim_{\mathbb{C}}\mathcal{S}_{1,\lambda}$ for every $\lambda{\in}\Par(r)$ and 
$\dim_{\mathbb{C}}\mathcal{S}_t=\dim_{\mathbb{C}}\mathcal{S}_1$, we have 
$\dim_{\mathbb{C}}\mathcal{S}_{t,\lambda}=
\dim_{\mathbb{C}}\mathcal{S}_{1,\lambda}$ for every $\lambda{\in}\Par(r)$. 
Thus we obtain that $\mathcal{S}_t=\bigoplus_{\lambda{\in}H(m,n;r)}\mathcal{S}_{t,\lambda}$. 
Since $\bar{\mathcal{A}}_q=(\mathcal{S}_q{\otimes_{R_0}}\mathbb{Q}(q))
{\otimes_{\mathbb{Q}(q)}}\overline{\mathbb{Q}(q)}$ and 
$\dim_{\mathbb{Q}(q)}\mathcal{A}_q=\dim_{\mathbb{C}}\mathcal{S}_t$, 
we immediately get 
$\dim_{\overline{\mathbb{Q}(q)}}\bar{\mathcal{A}}_q=\dim_{\mathbb{C}}\mathcal{S}_t$. 
Because $t$ is transcendental, 
$\dim_{\overline{\mathbb{Q}(q)}}\bar{\mathcal{A}}_{q,\lambda}
=\dim_{\mathbb{C}}\mathcal{S}_{t,\lambda}$ for every $\lambda{\in}\Par(r)$. 
Thus we conclude that 
$\bar{\mathcal{A}}_q=\bigoplus_{\lambda{\in}H(m,n;r)}\bar{\mathcal{A}}_{q,\lambda}$. 
The second statement is the direct consequence of the double centralizer theorem. 
\end{proof}

\end{document}